\documentclass[12 pt]{amsart}

\usepackage{amscd}
\usepackage[dvips]{graphicx}
\usepackage{amscd}
\usepackage{amsmath}
\usepackage{amsfonts}
\usepackage{amssymb}
\usepackage[latin1]{inputenc}

\theoremstyle{plain} \numberwithin{equation}{section}
\newtheorem{theorem}{Theorem}[section]

\newtheorem{lemma}[theorem]{Lemma}

\theoremstyle{definition}

\newtheorem{example}[theorem]{Example}

 \numberwithin{equation}{section}

\def\beqn{\begin{equation}}
\def\eeqn{\end{equation}}

\def\bE{\mathbf{E}}
\def\bV{\mathbf{V}}
\def\bc{\mathbf{c}}
\def\fa{\mathfrak a}

\def\fg{\mathfrak g}
\def\fk{\mathfrak k}
\def\fh{\mathfrak h}

\def\fp{\mathfrak p}
\def\fq{\mathfrak q}

\def\fz{\mathfrak z}

\def\Claminv2{|C(\Lambda)|^{-2}}

\def\bH{\mathbf{H}}

\newcommand{\gL}{\Lambda}
\newcommand{\gD}{\Delta}
\newcommand{\gO}{\Omega}

\newcommand{\R}{\mathbb{R}}
\newcommand{\C}{\mathbb{C}}

\newcommand{\N}{\mathbb{N}}
\newcommand{\Z}{\mathbb{Z}}
\newcommand{\Sy}{\mathrm{Sym}}

\newcommand{\Tr}{\mathop{\mathrm{Tr}} }
\newcommand{\Sym}{\mathop{\mathrm{Sym}} }

\newcommand{\SL}{\mathop\mathrm{SL}}

\newcommand{\cD}{\mathcal{D}}

\newcommand{\cL}{\mathcal{L}}

\newcommand{\GL}{\mathrm{GL}}
\newcommand{\Sp}{\mathrm{Sp}}

\newcommand{\SO}{\mathrm{SO}}

\newcommand{\su}{\mathfrak{su}}
\newcommand{\SU}{\mathrm{SU}}
\newcommand{\ad}{\mathrm{ad}}
\newcommand{\Ad}{\mathrm{Ad}}

\newcommand{\Ome}{\Omega}

\def\draft{\centerline{(Draft {\the \day}/{\the\month} \the \year.)}}
\newcommand{\set}[1]{\left\{#1\right\}}

\begin{document}
\title[The Generalized Segal-Bargmann transform]{The Generalized
Segal-Bargmann  transform and Special Functions}
\author{Mark Davidson and  Gestur \'{O}lafsson}
\address{Department of Mathematics, Louisiana State University, Baton Rouge, LA\ 70803,
USA} \email{davidson@math.lsu.edu} \email{olafsson@math.lsu.edu}
\thanks{Research by G. \'{O}lafsson supported by NSF grant DMS 0070607 and DMS 0139783}

\begin{abstract} Analysis of function spaces and special functions are
closely related to the representation theory of Lie groups. 
We explain here the connection between the Laguerre
functions, the Laguerre polynomials, and
the Meixner-Pollacyck polynomials on the one side,
and highest weight representations of Hermitian
Lie groups on the other side. The representation
theory is used to derive differential equations
and recursion relations satisfied
by those special functions.
\end{abstract}
\maketitle

\section*{Introduction}
\noindent The broad subject of special functions is
closely related to the representation theory of Lie groups and
arises very naturally in analysis, number theory, combinatorics,
and mathematical physics. The classical texts \cite{V68,VK91} by
Vilenkin and Vilenkin and Klimyk  well document this interplay. The interested
reader is also referred to  \cite{D80} by Dieudonn\'e and the
recent text \cite{AAR} by Andrews, Askey, and  Roy.    It
is not our aim to add to this discussion in a general way,  but
mainly to concentrate on our specific work \cite{df2,DO,doz1,doz2}
on the Meixner-Pollacyck polynomials, Laguerre functions, and Laguerre
polynomials defined on symmetric cones. It is in the interplay
among the various function spaces that model unitary highest
weight representation for a Hermitian group that new results
involving difference equations, differential equations, recursion
formulas, and generating functions arise.

The classical Laguerre polynomials $L_{n}^{\alpha}(x)$ can be
defined in several different ways. One of the oldest definitions
is in terms of the generating function
\begin{equation}\label{eq-generating}
(1-w)^{-\alpha-1}\exp\left(  \frac{xw}{w-1}\right)  =\sum_{n=0}^{\infty}%
L_{n}^{\alpha}(x)w^{n}\,,\quad |w|<1,\,-1<\alpha\,. %
\end{equation}

Another way is to use the Rodriguez type formula

\begin{equation}
L_{n}^{\alpha}(x)=\frac{e^{x}x^{-\alpha}}{n!}\frac{d^{n}\,\,}{dx^{n}}%
(e^{-x}x^{n+\alpha})\,.\label{e-rod}%
\end{equation}
Finally, they are special cases of hypergeometric functions:
\begin{eqnarray}
L_{n}^{\alpha}(x) &
=&\sum_{k=0}^{n}\frac{\Gamma(n+\alpha+1)}{\Gamma
(k+\alpha+1)}\frac{(-x)^{k}}{k!(n-k)!}\label{e-formula1}\label{eq1.3}\\
&  =&\frac{\Gamma(n+\alpha+1)}{\Gamma(n+1)}{}_{1}F_{1}(-n,\alpha
+1;x)\, .\label{e-formula}%
\end{eqnarray}

The polynomials
$$\sqrt{\frac{\Gamma (n+1)}{\Gamma (n+\alpha +1)}}\, L^\alpha_n(x),\quad n\in \mathbb{N},$$
form an orthonormal basis for the $L^2$-Hilbert space $L^2(\R^+,
e^{-x}x^{\alpha}dx)$, or otherwise stated:

\begin{theorem} The functions
\begin{equation}
\ell^\alpha_n (x)=L^\alpha_n(2x)e^{-x}
\end{equation}
form an orthonormal basis for $L^2(\R^+,d\mu_\alpha)$, where
$d\mu_\alpha (x)=x^\alpha\, dx$.
\end{theorem}

As was noted in \cite{doz1}, Theorem 3.1, a simple calculation
shows that the functions $$\set{\ell_n^\alpha,\; n\in
\mathbb{N}}$$ have Laplace transform
$$
\int_0^\infty e^{-xz}\ell^\alpha_n(x)\, d\mu_\alpha(x)=\frac{\Gamma (n+\alpha +1)}%
{\Gamma (n+1)}\left(\frac{z-1}{z+1}\right)^n(z+1)^{-(\alpha +1)}$$
which forms a basis for the  space of $\SO(2)$-finite vectors in a
holomorphic discrete series of a conjugate version of $\SL
(2,\R)$. These representation theoretic relationships
are the basic starting point in \cite{doz2} and will be described
in the following sections. In short let $\Omega\subset \R^d$ be a
symmetric cone and let $T(\Omega)=\Omega +i\R^d$ be the
corresponding tube domain. Let $G^c$ be a connected semisimple Lie
group, locally isomorphic to the connected component containing
the identity of the  group of holomorphic automorphisms of
$T(\Omega)$. We use a natural orthogonal set of $K$-finite vectors
in the Hilbert space, $\mathcal{H}_\nu(T(\Omega))$, of holomorphic
functions corresponding to a highest weight representation of
$G^c$ and the generalized Segal-Bargmann transform,
$\mathcal{L}_\nu$, defined in \cite{OOe96} to define an orthogonal
set $\{\ell_{m}^{\nu}(x)\}_m$ on the space
$L^{2}(T(\Omega),d\mu_{\nu})$ for an appropriate measure
$d\mu_{\nu}(x)$. A generalization of (\ref{eq1.3}) was used in
\cite{FK94}, p. 343, to define Laguerre polynomials,
$L_{m}^\nu$, for arbitrary symmetric cones and the
following relationship was obtained:
\begin{equation}
\ell_{m}^{\nu}(x)=e^{-\Tr(x)}L_{m}^{\nu}(2x)\, .\label{eq-ellL}%
\end{equation}
Applying the same ideas to the bounded realization of $T(\Omega)$
gives a natural generalization of the Meixner-Pollacyck
polynomials, \cite{doz2}, Section 5. Most of the main results in
\cite{DO, doz1,doz2}  concern the differential and
difference relations among these special functions. Several
authors have used the relation between highest weight modules and
the Meixner-Pollacyck polynomials/Laguerre polynomials. In
particular, we would like to mention
\cite{wG03,wGeK02,JJ98,eK01,KJ98,MR91,Z01,Z01a,Z02} as related
references, but as far as the authors know, the first time the
Laguerre polynomials show up as coefficient functions of a
representation of a Lie group was in \cite{V68}, p. 430--434,
where Vilenkin relates the Laguerre polynomials to coefficient
functions of representations of the group
$$\left\{\left(
\begin{matrix}1 & a & b\cr 0 & c & c\cr 0 & 0 &1\cr
\end{matrix}
\right)\mid a,b,c,d\in \mathbb{R}, c\not= 0\right\}\, . $$

This overview article is based on the presentation of the second
author at the conference \textit{Representations of Lie Groups,
harmonic analysis on homogeneous spaces and
quantization}, at the Lorentz Center, University of Leiden, December
9 -- 13, 2002. The second author would like to thank the
organizers G. van Dijk and V.F. Molchanov, as well as the
staff at the Lorentz Center for for the invitation and
hospitality during his stay in Leiden.

\section{Tube type domains}

\noindent In this section we introduce the necessary notation and
facts related to Hermitian Lie groups and bounded symmetric
domains $\cD$, containing the origin $0\in \C^d$. Even if much of
what we say is valid for all bounded symmetric domains in $\C^d$
we will assume from the beginning, as our applications will be
restricted to that case, that $\cD$ is isomorphic to a tube type
domain $T(\Omega)$. There exists a connected semisimple Lie group
that acts transitively on $\cD$. For simplicity we will assume
that $G$ is simple. The general situation can be reduced to this
case by considering direct products. Let $K=\{g\in G\mid g\cdot
0=0\}\, .$ Then $K$ is a maximal compact subgroup of $G$ and
$\cD\simeq G/K$. Furthermore there exists an involution $\theta
:G\to G$, a Cartan involution, such that
$$K=G^\theta =\{g\in G\mid \theta (g)=g\}\, .$$
In particular $\cD$ is a \textit{Riemannian symmetric space}. In
this case $\R^d$ has the structure of an Euclidean Jordan algebra
with identity element $e$, \cite{FK94}, and $\Omega =\{x^2\mid
x\in \R^d\}_0$ where the subscript ${}_0$ stands for the connected
components containing $e$. The isomorphism $\cD\simeq T(\Omega)$
is then given by the Cayley transform
$$\bc (z)=(e+z)(e-z)^{-1}=\frac{e+z}{e-z}.$$
Let $\fg_\C=\fg\otimes_\R\C$ denote the complexification of $\fg$
and denote by $G_\C$ a simply connected and connected complex Lie
group with Lie algebra $\fg_\C$. Then $g\cdot z$ is well defined
for $g\in G_\C$ and almost all $z\in \C^d$ and $(g_1g_2)\cdot
z=g_1\cdot (g_2\cdot z)$ whenever defined. In particular there
exists an element $c\in G_\C$ such that $c\cdot z=\bc (z)$ for all
$z\in \cD$. Assuming that $G\subset G_\C$ we set $G^c=cGc^{-1}$.
Then $G^c$ acts transitively on $T(\Omega)$ and the stabilizer of
$e\in \gO$ is given by $K^c=cKc^{-1}$. Finally, there exists a
conjugation $\sigma :\C^d\to \C^d$ which lifts to an involution
$\tau :G_\C\to G_\C$ such that, with $H=G^\tau_0$, we have
$$\cD_{\R} :=\{z\in \cD\mid \sigma (z)=z\}$$
is a totally real Riemannian subsymmetric space, isomorphic to
$H/H\cap K$. We can choose $\sigma$ so that the Cayley transform
maps $\cD_\R$ bijectively onto $\gO$. Notice that $\theta$ and
$\tau $ commute and hence
\begin{eqnarray*}
\fg&=&\fk\oplus \fp=\fh\oplus \fq\\
&=&\fh_k\oplus \fh_p\oplus \fq_k\oplus \fq_p
\end{eqnarray*}
where $\fp=\{X\in \fg\mid \theta (X)=-X\}$, $\fq=\{X\in \fg\mid
\tau(X)=-X\}$, and a subscript ${}_k$ (respectively ${}_p$)
indicates intersection with $\fk$ (respectively $\fp$). Everything
is now set up so that $\tau (g)= c^2(g)(c^{-1})^2$, $\fg^c=\Ad
(c)\fg$, and the Lie algebra of $G^c$ is given by
$$\fg^c=\fh_k\oplus i\fh_p\oplus i\fq_k\oplus \fq_p\, .$$
We note that $K^c:=cKc^{-1}=H_\C\cap G^c$ and
$H^c=cHc^{-1}=K_\C\cap G^c$.

The fact that $\cD$ is an irreducible bounded symmetric domains
implies that $\fz_\fk$, the center of $\fk$, is one dimensional.
We can choose $Z_0\in i\fz_\fk$ such that $\ad (Z_0):\fg_\C\to
\fg_\C$ has eigenvalues $0,1$ and $-1$. The $0$-eigenspace is
exactly $\fk_\C$. Define $\fp^+$ to be the $1$-eigenspace and
$\fp^-$ to be the $-1$-eigenspace.  Then $\fp^+$ and $\fp^-$ are
both abelian Lie algebras, normalized by $K_\C$, and
$\fg_\C=\fp^+\oplus \fk_\C\oplus \fp^-$. Denote by $P^+=\exp
(\fp^+)$ the closed Lie subgroup of $G_\C$ with Lie algebra
$\fp^+$ and similarly $P^-=\exp (\fp^-)$. Then $P^+\times
K_\C\times P^-\ni (p,k,q)\mapsto pkq\in G_\C$ is a diffeomorphism
onto an open dense subset of $G_\C$. Furthermore $G\subset P^+K_\C
P^-$. For $g\in P^+K_\C P^-$ we denote by $p(g)\in P^+$, $k_\C
(g)\in K_\C$, and $q(g)\in P^-$ the inverse of the above
diffeomorphism. Then the bounded realization of $G/K$ inside
$\fp^+$ is given by
$$G/K \ni gK\mapsto (\exp|_{\fp^{+}})^{-1} (p(g))\in \fp^+\, .$$

We will also need some basic information about roots. Let $\fa_k$
be a maximal abelian subspace of $\fq_k$. Then $\fa_k$ is in fact
maximal abelian in $\fq$. Let $\Delta$ be the set of roots of
$(\fa_k)_\C$ in $\fg_\C$ and let $\Delta_k$ be the set of roots of
$(\fa_k)_\C$ in $\fk_\C$. Then
$$\Delta_k=\{\alpha\in \Delta\mid \alpha(Z_0)=0\}$$
Set $\Delta_n=\{\alpha\in \Delta\mid \alpha (Z_0)\not=0\}$. Then
$\Delta_n=\Delta_n^+\cup \Delta_n^-$, disjoint union, where
$$\Delta_n^\pm =\{\alpha\in \Delta\mid \alpha (Z_0)=\pm 1\}\, . $$
Recall that two roots $\alpha,\beta\in \Delta$ are called
\textit{strongly orthogonal} if $\alpha\pm \beta\not\in \Delta$.
In our situation there are only two root lengths. Choose a maximal
set $\{\gamma_1,\ldots ,\gamma_r\}$ of long strongly orthogonal
roots in $\Delta_n^+$. Notice that $r=\dim \fa_k$. The above set
of roots can now be described as:
\begin{eqnarray}
\Delta_n^\pm &=&\pm [ \{\gamma_1,\ldots , \gamma_r\}\cup
\{\frac{1}{2}(\gamma_i+\gamma_j)\mid i< j\}]\label{eq-noncompact}\\
\Delta_k&=&\pm \{\frac{1}{2}(\gamma_i-\gamma_j)\mid i<
j\}]\label{eq-compact}
\end{eqnarray}
We choose the ordering such that $\gamma_1>\ldots >\gamma_r>0$.
Then the positive set of roots are those with a ` $+$' in the
above equations. For $\alpha \in \Delta$ let
$$s_\alpha : \gamma \mapsto \gamma -\frac{2(\gamma ,\alpha)}{(\alpha ,\alpha)}
\alpha$$ denote the corresponding Weyl group reflection. Here
$(\cdot ,\cdot )$ denotes a Weyl group invariant inner product on
$i\fa_k^*$ such that $(\gamma_1,\gamma_1)=1$. Notice that
$s_{(\gamma_i-\gamma_j)/2}(\gamma_i)=\gamma_j$ and
$s_{\gamma_i}((\gamma_i+\gamma_j)/2)=-(\gamma_i-\gamma_j)/2$. It
follows that all the root spaces $(\fg_\C)_{\gamma_j}$ have the
same dimension, which is in fact one, and also the spaces
$(\fg_\C)_{(\gamma_i\pm \gamma_j)/2}$ have the same dimension,
which we denote by $a$.

Let $H_j\in i\fa_k$ be such that $\gamma_i(H_j)=2\delta_{ij}$.
Then $i\fa_k=\oplus_{j=1}^r\R H_j$ and $Z=\frac{1}{2}(H_1+\ldots
+H_r)$. Let $\fa=\Ad (c)(\fa_k)=\Ad (c)^{-1}\fa_k$. Then $\fa$ is
maximal abelian in $\fh_p$ and $\fp$. The roots of $\fa$ in $\fg$
are given by $\Delta_\fa = \Delta\circ \Ad (c)$. We set:
\begin{eqnarray}
\label{eq-xi}\xi_j&=&\Ad(c)^{-1}H_j\cr
\xi& = &\xi_1+\ldots +\xi_r=2\Ad (c)^{-1}Z_0\nonumber\cr
\beta_j&=&\gamma_j\circ \Ad (c)\nonumber
\, .
\end{eqnarray}
Then $\fa=\bigoplus_{j=1}^r \R\xi_j$ and $\xi$ is central in $\fh$.
We use $\beta_1,\ldots ,\beta_r$ to identify $\fa_\C^*$ with
$\C^r$ by
$$(\lambda_1,\ldots ,\lambda_r)\mapsto \lambda_1\beta_1+
\ldots +\lambda_r\beta_r\, .$$
We also embed $\C$ into $\fa^*_\C$ by
$\lambda \mapsto \lambda (\beta_1+\ldots +\beta_r)$.

If $G$ is simple, then by the classification the triple $(G,K,H)$
is locally isomorphic to one of the following:
\bigskip

\begin{center}
\begin{tabular}[h]{|c|c|c|c|c|c|}
\hline $G$&$K$&$H$&$K\cap H$&$r$&$a$\cr\hline\hline
$\mathrm{Sp}(n,\R)$ &
$\mathrm{U}(n)$&$\mathrm{GL}(n,\R)_+$&$\mathrm{O}(n)$&$n$&$1$\cr\hline
$\mathrm{SU}(n,n)$&$\mathrm{U}(n)\times
\mathrm{U}(n)$&$\mathrm{GL}(n,\C)_+$&$\mathrm{U}(n)$&$n$&$2$\cr\hline
$\mathrm{SO}^*(4n)$&$\mathrm{U}(2n)$&$\mathrm{SU}^*(2n)\R^+$&$\mathrm{Sp}(2n)$&$2n$&$4$\cr\hline
$\mathrm{SO}(2,k)$&$\mathrm{S}(\mathrm{O}(2)\times
\mathrm{O}(k))$&
$\mathrm{SO}(1,k-1)\R^+$&$\mathrm{SO}(k-1)$&$2n$&$k-2$\cr\hline
$\mathrm{E}_{7(-25)}$&$\mathrm{E}_6\mathbb{T}$&$\mathrm{E}_{6(-26)}\R^+$&$\mathrm{F}_4$&$3$&$8$
\cr \hline\hline
\end{tabular}
\end{center}
\bigskip

Here the subscript ${}_+$ stands for the real and positive
determinant.

\begin{example}[Symmetric Matrices]
Let $\R^d=\mathrm{Sym}(n,\R)$, $d=n(n+1)/2$, be the real space of
$n\times n$ real symmetric matrices. The complexification of
$\R^d$ is then the space of $n\times n$ complex symmetric matrices
$\mathrm{Sym}(n,\C)$. Obviously $\mathrm{Sym}(n,\C)$ is a complex
Jordan algebra, where the product is given by $X\cdot
Y:=(XY+YX)/2$. Furthermore,
$$\gO=\{X\in \mathrm{Sym}(n,\R)\mid X>0\}\, ,$$
where $> $ stands for positive definite. The action of
$\mathrm{GL}(n,\R)$ on $\mathrm{Sym}(n,\R)$ is given by
$$g\cdot X= gXg^T\, .$$
The orbits are parametrized by the signature. In particular, the
cone of positive definite matrices is homogeneous. The identity
element is just the usual identity matrix $I_n\in \gO$ and the
stabilizer of $I_n$ is $\SO(n)$.

Write an element $g\in \mathrm{Sp}(n,\R)$ as
$$g=\left(\begin{matrix} A& B\cr C&D\end{matrix}\right)\, .$$
Then the action of $\mathrm{Sp}(n,\R)$ on $T(\Omega)$ is given by
\begin{equation}\label{eq-actSp}
g\cdot Z=(AZ+B)(CZ+D)^{-1}\, .
\end{equation}
Notice that (\ref{eq-actSp})  in fact defines an almost everywhere
defined action of $\Sp (n,\R)_\C=\Sp (n,\C)$ on $\Sy (n,\C)$. We
have
$$\cD=\{X\in \mathrm{Sym}(n,\C)\mid I_n-Z^*Z>0\}$$
and the Cayley transform is given by
$$c(Z)=(I_n+Z)(I_n-Z)^{-1}=2^{-2n}\left(\begin{matrix}I_n & I_n\cr
-I_n & I_n\end{matrix}\right)\cdot Z\, .$$ Define $\sigma
:\mathrm{Sym}(n,\C)\to \mathrm{Sym}(n,\C)$ to be the usual complex
conjugation. Then $\cD_\R=\{Z\in \Sy (n,\R)\mid I_n-Z^2>0\}$ and a
simple calculation shows that $\bc(X)\in\gO$ for $X\in \cD_\R$ as
was to be expected.

The case $\R^d=\{Z\in \mathrm{M} (n,\C)\mid Z^*=Z\}$, $G=\SU
(n,n)$ and $H=\GL (n,\C)_+$ is treated is the same way. Notice
that in this case the complexification of $\R^d$ is given by
$\R^d_\C=\mathrm{M}(n,\C)$ as every complex matrix can be written
in an unique way as $Z=X+iY$ with $X=X^*$ and $Y+Y^*$. Simply set
$X=1/2(Z+Z^*)$ and $Y=1/(2i)(Z-Z^*)$.
\end{example}

We refer to \cite{FK94,Loos-bsd,UP87} for further information on
the structure of bounded symmetric domains and the relation
between bounded symmetric domains and Jordan algebras. For the
more group theoretical description we refer to
\cite{He84,KW65a,KW65b}. The connection between
Hermitian symmetric spaces and causal symmetric
spaces $G/H$ is described in \cite{HO97}.

\section{Highest weight representations}
\noindent The unitary highest weight representations of $G$ are
well understood unitary representations that can be realized in a
Hilbert space of holomorphic functions on $\cD$. The idea
described in this section is to use the \textit{Restriction
principle}, introduced in \cite{OOe96} (see also
\cite{o00,oz-weyl}), to construct the Berezin transform and
generalized Segal-Bargmann transform and to transfer information
known from the highest weight representations to
extract information about
$L^2(\cD_\R)$ and $L^2(\Omega, d\mu_\nu)$, where $\mu_\nu$ is a
measure on $\Omega$ to be introduced in a moment. Thus: \textit{We want to
use the highest weight representations to do harmonic
analysis on bounded real domains and symmetric cones}. It should
be noted that almost all Riemannian symmetric spaces can be
realized as a real form of a Hermitian symmetric spaces. We will
in this section give a brief introduction to the theory of highest
weight representations. We refer to
\cite{EHW83,FK,HC56,HO92,HPJ83,N99,Ol81,RS86,RV76,W79} for more
information.

Let $p=2d/r$, and let $J(g,z)$ be the complex Jacobian determinant
of the action of $G$ on $\cD$. If we are discussing the tube type
realization, then we will use the same notation for the  complex
Jacobian determinant of the action of $G^c$ on $T(\Omega)$. The
map $(g,z)\mapsto J(g,z)$ can be expressed in terms of the
$K_\C$-projection $(g,z)\mapsto k_\C(g \exp(z))\in K_\C$ and
$$\rho_n=\frac{1}{2}\sum_{\alpha\in \gD_n} \dim
((\fg_\C)_\alpha)\alpha=\frac{1}{2}\left(1+\frac{a(r-1)}{2}\right)
(\gamma_1+\ldots +\gamma_r)$$
as
$$J(g,z)=\det \Ad (k_C(g\exp z))|_{\fp^+}=k_C(g\exp z)^{2\rho_n}\, .$$
Let $\gD$ be the determinant function on
the Jordan algebra, $\R^d$, and $\Tr$ the trace functional. In the
case $\R^d=\Sym (n,\R)$ this is just the usual determinant
function and trace. We denote by $\Delta_j$ the principal minors.
Recall that $\Delta_r=\Delta$.
For $(\alpha_1,\ldots ,\alpha_r)\in\C^r$ let
\begin{equation}\label{eq-Delta}
\Delta_\alpha(w):=\Delta_1(w)^{\alpha_1-\alpha_2}\gD_2(w)^{\alpha_2-\alpha_3}\ldots
\Delta_r(w)^{\alpha_r}
\end{equation} and
\begin{equation}\label{def-psi}
\psi_\alpha (x)=\int_{K\cap H}\Delta_\alpha (kx)\, dk\, .
\end{equation}
Let $h(z,w)=\gD (e-z\bar{w})$. Finally, we define the
Gindikin-Koecher gamma function by
\begin{equation}\label{eq-Gamma}
\Gamma_\Omega (\lambda ):=\int_\gO e^{-\Tr (x)}\Delta_\lambda
(x)\Delta^{-d/r}\, dx,\quad \lambda\in\C\, .
\end{equation}
We have
$$\Gamma_\Omega (\lambda )=(2\pi)^{(d-r)/2}\prod_{j=1}^r\Gamma (\lambda_j - \frac{a}{2}(j-1))\, ,$$
where $\Gamma$ is the usual $\Gamma$-function. For $m\in
\gL:=\{(m_1,\ldots ,m_r)\in \N_0^r\mid m_1\ge m_2\ge \ldots \ge
m_r\}$ denote by $\tau_{m}$ the irreducible $K\cap
H$-spherical representation with lowest weight $-m$.
Denote by $P(\C^r)$ the space of polynomial functions
on $\C^r$,
Then $\tau_m$ can be realized (with multiplicity one)
in a subspace $P_m(\C^r)\subset P(\C^r)$. Furthermore
(see \cite{Schmid})
$$P(\C^r)=\bigoplus_{m\in\gL}P_{m}(\C^r)\, .$$

There is a canonical series of highest weight representations (or
weighted Bargmann spaces in some cases)
$$(\pi_\nu,\bH_\nu(\cD))\quad \nu\in
\set{0,\frac a 2,\ldots, \frac a 2 (r-1)}\cup(\frac a
2(r-1),\infty),$$ the so-called Berezin-Wallach set.
To simplify the notation for the moment write
$\bH_\nu$ for $\bH_\nu(\cD)$. The space
$\bH_\nu$ is a Hilbert space of holomorphic functions on
$\mathcal{D}$. For  $\nu>a(r-1)+1$ the norm is given by
$$\|F\|_\nu^2=\alpha_\nu\int_{\mathcal{D} }^{ }|F(z)|^2 h(z,z)^{\nu-p}\,
dz,  \quad
\alpha_\nu=\frac{1}{\pi^d}\frac{\Gamma_\Omega(\nu)}{\Gamma_\Omega
(\nu - d/r)}\, .$$ In this case, if we assume that the center of
$G$ is finite, the representation $(\pi_\nu,\bH_\nu)$ can be
realized as a discrete summand in $L^2(G)$ (\cite{HC56}) and
$L^2(G/H)$ (\cite{OOe88,OOe91}). For $\nu<a(r-1)+1$ we use
analytic continuation of this norm.  The representation of $G$ on
$\bH_\nu$ is given by
\begin{equation}\label{eq-rep}
\pi_\nu(g)F(z)=J(g^{-1},z)^{\frac \nu p}F(g^{-1}z)\, .
\end{equation}

The facts that we will need are:
\begin{theorem} Let the notation be as above. Then the following holds:
\begin{enumerate}
\item $\bH_\nu$ is a reproducing kernel Hilbert space
with reproducing kernel
$$K_w(z)=K(z,w)=h(z,w)^{-\nu}=\Delta(e-z\bar{w})^{-\nu}.$$
\item If  $\nu>(r-1)\frac{a}{2}$ then the space of polynomials
$P(\C^d)$ is dense in $\bH_\nu$. More specifically, the space of
$K$-finite vectors $(\bH_\nu)_K$ can be naturally identified with
$P(\C^d)$
\item  All of the $K$-representations
$(\tau_{{m}},P_{{m}}(\C^d)$) are $K\cap H$-spherical and
  $$P_{ m}(\C^d)^{K\cap H}={\mathbb C}\psi_{m}.$$
  \item $\bH_\nu^{K\cap H}\simeq \bigoplus_{m\in\gL}{\mathbb
  C}\psi_{ m}.$
\item The norm of the function
$\psi_m$ is given by
$$\|\psi_{m} \|^2=d_m^{-1}\frac{\left(\frac{d}{r}\right)_m}{(\nu )_m}$$
where $d_m=\dim P_m(\C^d)$.
\end{enumerate}
\end{theorem}

Notice that (1) means that the point evaluations are continuous
linear functionals on $\bH_\nu$ and
$$F(z)=(F|K_z)\quad \text{for all } F\in
\bH_\nu.$$

\section{The Berezin transform and generalized Segal-Bargmann transform}
\noindent
In this section we recall the basic facts about the
\textit{Berezin transform} and \textit{generalized Segal-Bargmann
transform}. We refer to
\cite{B78,vanDijk-Pevsner,vDH97,vDM03,GM02,JO98,Neretin,o00,OOe96,oz-weyl,UU,UP83}
for further information. We define a map, \textit{the restriction
map}, (see \cite{OOe96}) $R_\nu:\bH_\nu\rightarrow
C^\infty(\mathcal{D}_{\mathbb R})$ by
$$R_\nu F(x)=h(x)^\frac{\nu}{2}F(x),$$ where $h(x)=h(x,x).$
Then $R_\nu$ is injective. Since the invariant measure $d\eta$ on
$\mathcal{D}_{\mathbb R}$ is $d\eta(x)=h(x)^{-\frac {p}{2}}\; dx$
(and a few other things) we get
\begin{lemma}
$R_\nu F \in L^2(\mathcal{D}_{\mathbb R},d\eta)$ for all $F\in
P(\C^d)$ if and only if $\nu> \frac a 2 (r-1)$.
\end{lemma}

Thus the restriction map $$R_\nu: (\bH_\nu)_K =P(\C^d)\rightarrow
L^2(\mathcal{D}_{\mathbb R},d\eta)$$ is defined for all $\nu$ in
the continuous part of the Berezin-Wallach set.  We collect
the important properties in the following theorem (see
\cite{doz2}, section 3):
\begin{theorem}
Assume that $\nu> \frac a 2 (r-1)$.  Then the following hold:
\begin{enumerate}
  \item $R_\nu(P(\C^d))$ is dense in $L^2(\mathcal{D}_{\mathbb
  R},d\eta).$
  \item $R_\nu:\bH_\nu\rightarrow L^2(\mathcal{D}_{\mathbb
  R},d\eta)$ is closed.
  \item  The Berezin transform $R_\nu R_\nu^*$ is given by
  $$R_\nu R_\nu^*f(y)=\int_{\mathcal{D}_{\mathbb R}}
\frac{h(y)^{\frac \nu 2}h(x)^{\frac \nu 2}}{h(y,x)^\nu}
  f(x)\, d\eta(x)=D_\nu\star f(h)\quad (y=h\cdot 0),$$ where
  $D_\nu(h)=J(h,0)^{\frac {\nu}{p}}.$
  \item If $\nu>a(r-1)$ then $R_\nu R_\nu^*:L^2(\mathcal{D}_{\mathbb
  R},d\eta)\rightarrow L^2(\mathcal{D}_{\mathbb R},d\eta)$ is continuous
  with norm $||R_n R_n^*||_2\le ||D_\nu||_{L^1}<\infty.$
  \item If $\nu>a(r-1)$ then $R_\nu R_\nu^*:L^\infty(\mathcal{D}_{\mathbb
  R},d\eta)\rightarrow L^\infty(\mathcal{D}_{\mathbb R},d\eta)$ is continuous
  with norm $||R_n R_n^*||_\infty\le ||D_\nu||_{L^1}<\infty.$
\end{enumerate}
\end{theorem}
We are thus able to apply $R_\nu R_\nu^*$ to bounded spherical
functions!  The proof of (3) is standard:
\begin{eqnarray*}
R_\nu^*f(z)&=&(R_v^*f|K_z)\\
 &=&(f|R_vK_z)\\
  &=&\int_{ }^{ }f(x) h(x)^{\frac \nu 2} h(y,x)^{-\nu}\, d\eta.
\end{eqnarray*}
This implies $$R_\nu R_\nu^* f(y)=\int_{ }^{ }f(x)h(y)^{\frac \nu
2} h(x)^{\frac \nu 2} h(y,z)^{-\nu}\, d\eta.$$ In particular, this
shows that the restriction principle defined in \cite{OOe96}
results in the Berezin transform.

For $\nu>\frac  a 2 (r-1)$, the restriction map $R_\nu: \bH_\nu
\rightarrow L^2(\mathcal{D}_{\mathbb R},d\eta)$ is closed with
dense image. We can therefore define a unitary isomorphism $U_\nu:
L^2(\mathcal{D}_{\mathbb R},d\eta)\rightarrow \bH_\nu$ by
polarization. Thus $U_\nu$ satisfies
$$R_\nu^*=U_\nu\sqrt{R_\nu R_\nu^*}.$$
The map $U_\nu$ is the generalized Segal-Bargmann transform. Let
$c(\nu)$ denote the \textit{Harish-Chandra} $c$-function, $W$ the
Weyl group corresponding to the root system $\Delta(\fa,\fh)$,
$w=\# W$, and  $\mathcal{F }$ the Harish -Chandra/ Helgason
spherical Fourier transform
\begin{eqnarray*}\label{def-Fourier}
\mathcal{F}:L^2(\mathcal{D}_{\mathbb R},d\eta)^{H\cap
K}&\rightarrow& L^2(\mathfrak{a},\frac 1 \omega
\frac {d\lambda}{|c(\lambda)|^2})^W\\
F&\mapsto& \int_{\cD_\R}F(x)\varphi_\lambda (x)d\eta (x)
\end{eqnarray*}
where $\varphi_\lambda$ is the spherical functions on $H/H\cap
K\simeq \cD_\R$ given by
$$\varphi_\lambda (x)=\psi_{i\lambda +\rho}((e+x)(e-x)^{-1})=\psi_{i\lambda +\rho}(\bc (x))\, .$$
Here $\rho$ is half the sum of the positive roots for $H/(H\cap
K).$ Then the above discussion results in a unitary isomorphism
$$\mathcal{F}\circ U_\nu^*: P(\C^r)^{K\cap H}=\bH_\nu^{H\cap K}\rightarrow
L^2(\mathfrak{a},\frac 1 w \frac{d\lambda}{|c(\lambda)|^2})^W$$
and we get:

\begin{theorem}
If $\nu>\frac a 2 (r-1)$ then the functions
$$\set{\mathcal{F}\circ U_\nu^*(\psi_{ m })}_{{\mathbf
m}\ge 0}$$ form an orthogonal basis for $L^2(\mathfrak{a},\frac 1
w \frac{d\lambda}{|c(\lambda)|^2})^W.$
\end{theorem}

\section{The Polynomials $p_{\nu,{ m}}(\lambda)$}
\noindent The next obvious task is to understand the functions
$\mathcal{F}\circ U_\nu^*(\psi_{ \mathbf m })$.  For that we
define the polynomials $p_{\nu,{ \mathbf m}}(\lambda)$ by
\begin{equation}\label{de-pnu}
p_{\nu,m}(\lambda)=\|\psi_{\mathbf{
m}}\|^{-2}\psi_{m}(\partial_x)[\Delta(e-x^2)^{-\frac \nu
2}\phi_\lambda(x)]_{x=0},
\end{equation}

We  have:

\begin{lemma}[\cite{doz2}, Lemma 4.1]
$$\Delta(e-x^2)^{-\frac \nu 2}\phi_\lambda(x)=\sum_{m }^{ }p_{\nu,m}(\lambda)\psi_{ m}(x).$$
\end{lemma}

We define the constants $c_\nu$ and $b_\nu$ by
$$c_\nu=\|D_\nu\|_{L^1}\quad \mbox{ and }\quad R_\nu R_\nu^*(\varphi_\lambda)
=c_\nu^{-1}b_\nu (\lambda )\varphi_\lambda \, .$$ The numbers
$c_\nu$ and $b_\nu (\lambda)$ have been evaluated by Zhang,
\cite{Z01a}.

\begin{theorem}[\cite{doz2}, Proposition 4.3]
Assume $\nu > \frac a 2 (r-1)$. Then $$ \mathcal{F}(U_\nu^* \psi_{
\mathbf m})(\lambda)=c_\nu^{-\frac 1 2}\sqrt{b_\nu(\lambda)}
\|\psi_{ \mathbf m}\|_\nu^2 p_{\nu,{\mathbf m}}(\lambda)\, .$$
\end{theorem}
The numbers $c_\nu$ and $b_\nu(\lambda)$ have been computed by
Zhang.

\section{Example: $G=SU(1,1)$}
\noindent Let $\mathcal{D}=\set{z\in {\mathbb C}\mid |z|<1}$.
Whenever defined we write for $g=\left(\begin{matrix}a & b\cr c&d
\end{matrix}\right)\in \GL (2,\C)$ and $z\in \cD$:
$$g\cdot z =\frac{a z +b}
{c z +d}\, .$$ Then, restricted to $\SU (1,1)$, this defines a
transitive action of $\SU (1,1)$ on $\cD$. In this case we take
the conjugation $\sigma$ as the usual complex conjugation
$z\mapsto \bar{z}$. Then $\mathcal{D}_{\mathbb R}=(-1,1)$ and
$\psi_m(x)=x^m$. Notice that
$$H=\left\{h_t=\left(\begin{matrix}
\cosh (t) &\sinh(t)\cr \sinh (t) & \cosh(t)\end{matrix} \right)\mid
t\in \R\right\}$$ and $h_t(0)=\tanh (t)$. In particular $K\cap
H=\{\mathrm{Id}\}$ is trivial. As $(1+x)(1-x)^{-1}$ is real and
positive for all $x\in (-1,1)$ we can define
$$G_{\nu,\lambda}(x)=(1-x^2)^{-\frac \nu 2}(\frac {1 +
x}{1-x})^{i\lambda}\, .$$ Then we can expand $G_{\nu,\lambda}$ as
$$G_{\nu,\lambda}(x)=\sum_{n=0}^{\infty}p_{n,\nu}(\lambda)x^n,$$
with $p_{n,\nu}(\lambda)=(\frac \nu 2 + i\lambda)_n \;
{}_2F_1(\frac \nu 2 -i\lambda, -n, -\frac \nu 2 - i\lambda -n +1 ,
1 ),$ the Meixner-Pollacyck polynomials.

The Hilbert space $\bH_\nu$ is given as the space of holomorphic
functions $F:\cD\to \C$ such that
$$\|F\|^2=\frac{1}{\pi}\frac{\Gamma (\nu)}{\Gamma (\nu -1)}\int |F(x+iy)|^2(1-(x^2+y^2))^{2\nu-2} dxdy$$
and, whenever defined for $g=\left(\begin{matrix}a & b \cr c &
d\end{matrix}\right)$ we have
$$\pi_\nu (g)F(z)=
(-bz+a)^{-2\nu}f(g\cdot z)$$

The connection to the representation $\pi_\nu$ is as follows.  We
have $\xi=\begin{pmatrix}
  0 & 1 \\
  1 & 0
\end{pmatrix}\in \mathfrak{h}$ and  $Z_0=\frac 1 2\begin{pmatrix}
  1 & 0 \\
  0 & -1
\end{pmatrix}\in \mathfrak{k}$. Thus
\begin{eqnarray*}
\pi_\nu(\xi )f(x)&=&\nu x f(x)- (1-x^2)f'(x)\, , \cr
\pi_\nu(\xi)z^m&=&(\nu +m)z^{m+1}-mz^{m-1}\, ,\cr
\pi_\nu(\xi)G_{\nu,\lambda}&=&-2i\lambda G_{\nu,\lambda}\, .
\end{eqnarray*}

The last item follows from the simple calculation
\begin{equation}\label{surecurrence}
2i\lambda
p_{\nu,n}(\lambda)=(n+1)p_{\nu,n+1}(\lambda)-(\nu+n-1)p_{\nu,n-1}(\lambda)\,
.
\end{equation}

Now using the following facts that
\begin{eqnarray*}
\pi_\nu(-Z_0)f(x)&=& xf'(x)\, ,\cr \pi_\nu(-Z_0)z^m &=& mx^m\,
,\cr xG_{\nu,,\lambda}'(x)&=&\sum_{m=0}^{\infty}m
p_{m,\nu}(x)x^m\, .
\end{eqnarray*}
we get
\begin{equation}\label{eqgnu}
-(\nu+2n)p_{\nu,n}(\lambda)=(-\frac{\nu}{2}
+i\lambda)p_{\nu,n}(\lambda+i)-(\frac{\nu}{ 2 }+ i
\lambda)p_{\nu,n}(\lambda-i)\, .
\end{equation}

\section{Difference Relations}
\noindent The general recurrence relation corresponding to
(\ref{surecurrence}) and difference relation generalizing
(\ref{eqgnu}) are as follows:  Choose $Z_0\in
\mathfrak{z}_{\mathfrak{k}_{\mathbb C}}$ as before, i.e.,  such
that $\mathrm{ad}(Z_0)$ has eigenvalues $0, 1$, and  $-1$.  As
before set $\xi=\mathrm{Ad}(c)(-2Z_0)=\Ad (c)^{-1}(2Z_0)$. Denote
by $e_k$ the standard basis vector of $\R^r$ with a $1$ in each
$k$-th position and $0$'s elsewhere. For ${n}$ an $r$-tuple
define
$$\left(\begin{matrix} {n}\cr
{n}-e_k\end{matrix}\right):=(n_k+\frac{a}{2}(r-k))\prod_{j\not=
k}\frac{n_k-n_j+\frac{a}{2}(j-k-1)}{ n_k-n_j+\frac{a}{2}(j-k)}$$
and
$$c_{n} (k)=\prod_{j\not= k}\frac{n_k-n_j-\frac{a}{2}(j-k-1)}{n_k-n_j-\frac{a}{2}(j-k)}\, .$$
Then

\begin{theorem}[\cite{doz2},Theorem 5.2]
With the above notation we have:
$$\pi_\nu(-\xi)\psi_{ m}
=\sum_{j=1}^{r}\left(\begin{matrix} m\cr
m-e_j\end{matrix}\right)\psi_{{{ m}}-e _j}
-\sum_{j=1}^{r}(\nu +m_j -\frac a 2 (j-1))c_{ m}(j)\psi
_{{ m}+e_j}.$$
\end{theorem}

Let $q_{m,\nu}(z)=\Delta(z+e)^{-\nu}
\psi_{m}\left( (z-e)(z+e)^{-1}\right)$. Since $\xi =-2\Ad
(Z_0)$ and $Z_0$ is central in $\fk$ it follows, by some
calculations, that:
\begin{eqnarray*}
\pi_{\nu}(\xi)q_{m,\nu}&=&(r\nu+2\|m|)q_{m,\nu}\,
\cr
\pi_{\nu}(-2Z_{0})q_{m,\nu}&=&\sum_{j=1}^{r}\left(\begin{matrix}m\cr
m-e_{j}
\end{matrix}\right)\, q_{m-e_{j},\nu}-\sum_{j=1}^{r}%
(\nu+n_{j}-\frac{a}{2}(j-1)) c_{m}(j)q_{m+
e_{j},\nu}.
\end{eqnarray*}
Putting those pieces together we get:

\begin{theorem}[\cite{doz2}, Theorem 5.6]
We have the following difference relations amongst the
$p_{\nu,m}(\lambda)$:

$$2\sum_{j=1}^{r}(i\lambda_j + \rho_j)p_{\nu,m}(\lambda)=\sum_{j=1}^{r}\begin{pmatrix}
  m+e_j \\
  m
\end{pmatrix} p_{\nu,m+e_j}(\lambda)-(\nu +m_j -1
-\frac{a}{2}(j-1))c_{m-e_j}p_{\nu,
m-e_j}(\lambda)\, .$$
\end{theorem}
\begin{theorem}[\cite{doz2}, Theorem 6.1]
We have
\begin{eqnarray*}
-(r\nu+2|
m|)p_{\nu,m}(\lambda)&=&\sum_{j=1}^{r}\begin{pmatrix}
  i \lambda+\rho-\frac \nu 2 \\
  i \lambda + \rho -\nu 2 -e_j
\end{pmatrix}p_{\nu, m}-\cr
&&\quad\quad \sum_{j=1}^{r}(\frac{\nu}{2}+ i \lambda_j +\rho_j
-\frac{a}{ 2} (j-1))c_{i \lambda +\rho -\frac \nu 2}(j)p_{\nu,
m}(\lambda-ie_j).
\end{eqnarray*}
\end{theorem}

\section{The Unbounded Realization}
\noindent In this section we discuss the unbounded realization
$\mathcal{D}\simeq T(\Omega)= i{\mathbb R}^d+\Omega$. This can be
used to study recurrence and differential equations for Laguerre
polynomials and functions on $\Omega$.  We refer to
\cite{doz2,FK94,RV} for further references.

From now on $\nu\in \C$ is identified with $\nu (\beta_1+\ldots +\beta_r)\in
\fa^*$.
Using the Cayley transform we get a space of holomorphic functions
on $T(\Omega)$, ${\bH}_\nu(T(\Omega))=\pi_\nu(c)\bH_\nu$, and get
an orthonormal basis $$q_{{
m},\nu}(z)=\Delta(z+e)^{-\nu}\psi_{ m}(\frac{z-e}{z+e})$$ for
$\widetilde{\bH}_\nu^{H\cap K}$.  Notice that these functions
correspond to the functions $z\rightarrow
\frac{(z-i)^m}{(z+i)^{\nu+m}}$ on the upper half plane in the case
$G=\SL (2,{\mathbb R})$. Let us describe this in  more detail, as
the standard normalization for the unbounded realization is
usually different from the one in the bounded realization. For
$\nu>1+a(r-1)$ let $\bH_{\nu}(T(\Omega))$ be the space of
holomorphic functions $F:T(\Omega)\rightarrow\mathbb{C}$ such that
\begin{equation}
  \|  F\|    _{\nu}^{2}:=\beta_{\nu}\int_{T(\Omega
)}|F(x+iy)|^{2}\Delta(y)^{\nu-2d/r}\,dxdy<\infty\label{eq-normub}%
\end{equation}
where
\begin{equation}
\beta_{\nu}=\frac{2^{r\nu}}{(4\pi)^{d}\Gamma_{\Omega}(\nu-d/r)}\,.
\label{eq-cnu}%
\end{equation}
Then $\bH_{\nu}(T(\Omega))$ is a non-trivial Hilbert space. For
$\nu\leq 1+a(r-1)$ this space reduces to $\{0\}$ and analytic continuation
is again used to define the norm \cite{RV}. If $\nu=2d/r$
this is the \textit{Bergman space}.

Instead of using the $H^c$-invariant measure on the cone $\gO$ we
use a weighted measure $d\mu_\nu (x)=\gD (x)^{\nu -d/r}dx$
corresponding to the measure $x^{\nu -1}dx$ on $\R^+$. Thus we get
the weighted $L^2$-spaces
$$L^2_\nu(\Omega)=L^2(\Omega,\Delta(x)^{\nu-\frac d r}\; dx).$$
In this normalization the restriction map is simply $R_\nu
F=F|_\Omega$ and the unitary part becomes the Laplace transform
$\cL_\nu$.
We define $\mathcal L_\nu$ on the domain
$$
\{f\in L^2_\nu(\gO,d\mu_\nu)\mid  \int_\Omega e^{-(\omega ,x)}
|f(x)|d\mu_{\nu}(x)<\infty \}
$$
for all $\omega\in \gO$ by
$$\cL_\nu(f)(z):= \int_\Omega e^{-(z,x)} f(x)d\mu_\nu (x)\, .$$
Note that by the Cauchy-Schwarz inequality, the condition $f\in
L^2_\nu(\Omega)$ implies,  since $|e^{-(s+it, x)}|=e^{-(s,x)}$,
that $\mathcal L_\nu f$ is a well-defined
function on $T(\Omega)$. To simplify notation we
sometimes write $R$ for $R_\nu$ and $\cL$ for $\cL_\nu$.
Denote by $\cL^\gO=R_\nu\circ \cL_\nu$. Then $\cL^\gO$ is a
self-adjoint positive operator $L^2(\gO,d\mu_\nu)\to
L^2(\gO,d\mu_\nu)$.

\begin{theorem}\label{th-hwunb}
Let the notation be as above. Assume that $\nu>1+a(r-1)$. Then the
following hold:

\begin{enumerate}
\item  The space $\bH_{\nu}(T(\Omega))$ is a reproducing kernel Hilbert space.
\item The map
$$\Psi_\nu :=\frac{1}{\sqrt{\Gamma_\gO (\nu )}}\, \pi_\nu (c)^{-1}: \bH_\nu(T(\Omega))
\to \bH_\nu(\cD)$$ is a unitary isomorphism.

\item  The reproducing kernel of $\bH_{\nu}(T(\Omega))$ is given by
\[
K_{\nu}(z,w)=\Gamma_{\Omega}(\nu)\Delta\left(  z+\bar{w}\right)  ^{-\nu}%
\]

\item  If $\nu>(r-1)\frac{a}{2}$ then there exists a Hilbert space
$\bH_{\nu}(T(\Omega))$ of holomorphic functions on $T(\Omega)$ \
such that $K_{\nu}(z,w)$ defined in (2) is the reproducing kernel
of that Hilbert space and the universal covering  group of $G^{c}$
acts
unitarily and irreducibly on $\bH_{\nu}%
(T(\Omega))$.

\item  The map
\[
L_{\nu}^{2}(\Omega)\ni f\mapsto F=\mathcal{L}_{\nu}(f)\in\bH_{\nu
}(T(\Omega))
\]
is a unitary isomorphism and

\item  If $\nu>(r-1)\frac{a}{2}$ then the functions
\[
q_{m,\nu}(z):=\Delta(z+e)^{-\nu}\psi_{m}\left(
\frac
{z-e}{z+e}\right),  \,\qquad m\in{\Lambda},%
\]
form an orthogonal basis of $\bH_{\nu}(T(\Omega))^{L}$, the space
of $L$-invariant functions in $\bH_{\nu}(T(\Omega))$.
\end{enumerate}
\end{theorem}

The proof of (5) follows from a restriction principle argument.
Let $f\in L^2_\nu(\Omega)$ be in the domain of $RR^{\ast}$, then
\begin{equation}
\label{RR-ast}
\begin{split}
RR^{\ast}f(y)
&  =R^{\ast}f(y)\\
&  =(R^{\ast}f\mid K_{y})_{\bH_{\nu}(T(\Omega))}\\
&  =(f, RK_{y})_{L_{\nu}^{2}(\Omega)}\\
&  =\int_{\Omega}f(x)\overline{K(x,y)}\,d\mu_{\nu}(x)\\
&  =\Gamma_{\Omega}(\nu)\int_{\Omega}f(x)\Delta(x+y)^{-\nu}\,d\mu_{\nu}(x)\\
&  =\int_{\Omega}f(x)\mathcal{L}^\Omega (e^{-(y, \cdot)})(x)\,d\mu_{\nu}(x)\\
&  =\int_{\Omega}e^{-(y, x)}\Delta(x)^{\nu-d/r}\mathcal{L}^\Omega(f)(x)\,dx\\
&={\mathcal L^\Omega} (\mathcal L^\Omega f)(y),
\end{split}
\end{equation}
and
$$(\mathcal L^\Omega f,
\mathcal L^\Omega f)=(R^{\ast}f, R^{\ast}f)<\infty.
$$
Thus  $f$ is in the domain of $(\mathcal L^\Ome)^2$ and
$RR^{\ast}=(\mathcal L^\Ome)^2$. Therefore $(\mathcal L^\Ome)^2$
is a self-adjoint extension of $RR^\ast$, which is also
self-adjoint by the von Neumann theorem.

Consider the inverse operator $R^{-1}$ acting on the image of $R$.
For a function $g$ in the image of $R$,  $R^{-1} g$ is the unique
extension of $g$ to a holomorphic function on $T(\Ome)$. Thus
$R^{-1}{\mathcal L^\Ome}={\mathcal L}_\nu$. Letting $R^{-1}$ act
on the previous equality (\ref{RR-ast}) we get
$$
R^\ast f={\mathcal L}_\nu {\mathcal L^\Ome} f.
$$
This proves the polar decomposition formula. Since $R^\ast $ is
densely defined and $R$ is an injective closed operator we have
that the unitary part ${\mathcal L}_\nu$ extends to a unitary
operator. Thus we get the well known fact, first proved by Rossi
and Vergne \cite{RV} for the full Wallach set:
\begin{theorem}\label{th-laplace}
Let $\nu>a(r-1)/2$. Then the Laplace transform $\cL :
L^2(\Omega,d\mu_\nu) \to \bH_\nu (T(\Omega)$
$$\cL (f)(z)=\int_{\Omega}e^{-(z,x)}f(x)\Delta (x)^{\nu-d/r}\, dx$$
is an unitary isomorphism.
\end{theorem}

Combinging Theorem \ref{th-laplace} and Theorem  \ref{th-hwunb},
part 2, we get the following simple, but usefull fact:

\begin{lemma}\label{le-uniiso}
The map $\Xi=\Xi_\nu : L^2(\gO,d\mu_\nu)\to \bH(\cD)_\nu$,
\begin{eqnarray*}
\Xi (f)(w)&=&\frac{1}{\sqrt{\Gamma_\gO (\nu)}}\, \pi_\nu (c)^{-1}\cL (f)(w)\\
&=&\sqrt{\frac{2^{r\nu}}{\Gamma_\gO (\nu)}}\, \int_\gO e^{-(c
w,x)}f(x)\Delta(x)^{\nu-d/r}dx
\end{eqnarray*}
is a unitary isomorphism.
\end{lemma}

\section{The Laguerre Polynomials and Functions}
\noindent
In this section $\nu$ will stand for a complex number
identified with the element $\nu(\beta_1+\ldots +\beta_r)\in \fa^*$.
We define, as in \cite{FK94} (c.f. equation
(\ref{eq1.3}) in the Introduction), the Laguerre polynomials by
$$L_{m}^\nu(x)=(\nu)_{m} \sum_{|{{n}}|\le |{m}|}^{
}\begin{pmatrix}
{m}  \\
   {{ n}}
\end{pmatrix} \frac {1}{(\nu)_{{n}}}\psi_{{n}} (-x)$$
and the Laguerre functions $$\ell_{
m}^\nu(x)=e^{-\mathrm{Tr}(x)}L_{ m}^\nu(2x)$$
where
$$(\nu)_m=\frac{\Gamma_\Omega (\lambda +m)}{\Gamma_\Omega (\lambda)}\, .$$
Then
\begin{theorem}[\cite{doz2} Theorem 7.8]
$$\mathcal{L}_\nu(\ell_{ m}^\nu)=\Gamma_\Omega({
m}+\nu)q_{{{m}},\nu}.$$ In particular, it follows
that $\set{\ell_{{ m}}^\nu}_{{ m}} $ is an
orthogonal basis for $L^2(\Omega,d\mu_\nu)^{H\cap K}$.
\end{theorem}
It follows from the fact
$\|\psi_m\|^2=\frac{\left(\frac{d}{r}\right)_m}{d_m(\nu
)_m}$ and Lemma \ref{le-uniiso}, that
\begin{lemma}\label{le-norm}
$$\|\ell^\nu_m\|^2=\frac{\Gamma_\Omega (\nu )\left(\frac{d}{r}\right)_m}%
{2^{r\nu}d_m(\nu )_m}.$$
\end{lemma}

Furthermore the action $\pi_\nu(2Z_0)$ and $\pi_\nu(\zeta)$
translates into:
\begin{theorem}[\cite{doz2}, Theorem 7.9]
\begin{enumerate}
\item[]
\item $-(\nu r +2E)\ell_{{ m}}^\nu=\sum_{j=1}^{r} \begin{pmatrix}
{{ m}}  \\
{{ m}}-e_j \
\end{pmatrix}  ({m}_j -1 +\nu -\frac a 2 (j-1))
\ell_{{m}- e_j}^\nu +
\sum_{j=1}^{r}c_{m}(j)\ell_{{m}+e_j}^\nu$

\item $\pi_\nu(\zeta)\ell_{m}^\nu=(r \nu +2 |{m}|)\ell_{m}^\nu.$
\end{enumerate}
\end{theorem}

There are some things that remain to be done.  First, one can show
that $\mathfrak{g}_{\mathbb C}\simeq \mathfrak{sl}(2,{\mathbb
C})\simeq {\mathbb C}X^-\oplus {\mathbb C}X^+\oplus {\mathbb
C}Z_0$, where $X^-\in \mathfrak{p}^+$ corresponds to an
annihilating operator, $X^+\in \mathfrak{p}^-$ corresponds to  a
creation operator, and $Z_0\in \mathfrak{k}_{\mathbb C} $
corresponds to the Laguerre operator, when viewed as operators on
$L^2(\Omega,d\mu_\nu)$.  The problem is to find explicit formula
as second order differential operators and the full representation
of $\mathfrak{g}_{\mathbb C}$ on  $L^2(\Omega,d\mu_\nu)$. This has
been worked out for $G=SU(n,n)$\cite{DO}.

\begin{theorem}[\cite{DO}, Theorem 6.1]
For the cone of positive definite Hermitian matrices we have
\begin{enumerate}
  \item $\mathrm{tr}(-s\nabla\nabla -v\nabla +s)\ell_{ m}^\nu=(r\nu
  +|{{ m}}|)\ell_{ m}^\nu$
  \item $ \frac 1 2 \mathrm{tr} (s \nabla \nabla +(\nu I +2s)\nabla +(v I
  +s)\ell_{m}^\nu = \sum_{j=1}^{r} \begin{pmatrix}
     {{ m}}  \\
     {{ m}}-e_j \
  \end{pmatrix}(m_j-1 +\nu -\frac a 2 (j-1))\ell_{{{ m}}-e_j}^\nu$
  \item  $\frac 1 2 \mathrm{tr} (-s \nabla \nabla +(-\nu I +2s)\nabla +(v I
  -s))\ell_{{ m}}^\nu = \sum_{j=1}^{r}c_{ m}(j)\ell_{{
  m}+e_j}^\nu$.
\end{enumerate}
\end{theorem}

These formulas generalize the classical relations
\begin{enumerate}
  \item $(tD^2 +\nu D-t)\ell_n^v=-(2n+\nu)\ell_n^\nu$
  \item $(tD^2 + (2t+\nu)D
  +(t+\nu))\ell_n^\nu=-2(n+\nu-1)\ell_{n-1}^\nu$
  \item $(tD^2 - (2t-\nu)D
  +(t-\nu))\ell_n^\nu=-2(n+1)\ell_{n+1}^\nu$.
\end{enumerate}

\section{The generating function for the Laguerre polynomials}
\noindent Recall equation (\ref{eq-generating})
\begin{equation}\label{eq-generating1}
(1-w)^{-\alpha-1}\exp\left(  \frac{xw}{w-1}\right)  =\sum_{n=0}^{\infty}%
L_{n}^{\alpha}(x)w^{n}\,,\quad |w|<1,\,-1<\alpha\,. %
\end{equation}
from the introduction. Let us rewrite this using our notation $\nu
=\alpha+1$ and $\ell^\alpha_m (x)=L^\alpha_m(2x)e^{-x}$. By
replacing in (\ref{eq-generating1}) $x$ by $2x$ and multiplying
both sides by $\exp(-x)$ we obtain
\begin{equation}\label{eq-generating2}
(1-w)^{\nu}\exp\left(x\frac{1+w}{1-w}\right) =\sum_{m=0}^\infty
\ell^\nu_m (x)w^m
\end{equation}
We will now show how to generalize (\ref{eq-generating2}) using
the connection with the highest weight modules.

Suppose $\mathbb{H}$ is an arbitrary highest weight space for a
unitary representation of $G$ isomorphic to $\bH_\nu$. Let
$\bV=\mathbb{H}^{\fp^+}$ be
the lowest $K$-type. In other words
$$\bV=\{v\in
\mathbb{H}\mid X\cdot v=0, \text{for all $X\in \mathfrak{p}^+$}\}\, .$$
In this case $\bV=\mathbb{C}v_\circ$ is one dimensional and we can
thus identify $H$ and $\mathbb{C}$. Assume $v_\circ$ has norm $1$.
For each $T\in \mathfrak{p}^+$ we define a map $q_T:H\rightarrow
\mathbb{H}$ by the formula
$$q_T(v)=\sum_{m=0}^\infty \frac {\overline{T^n}\cdot v}{n!}.$$
We have the following
\begin{theorem}[\cite{df} Theorems 5.1, 7.1 and 7.2]
\begin{enumerate}\item[]
\item The series that defines $q_T$ converges in $\mathbb{H}$ if
and only if $T\in \mathcal{D}$.
\item  The map $\Xi:\mathbb{H}\rightarrow \bH_\nu$ given by
$$\Xi\xi(T)=q_T^*\xi,\quad T\in \mathcal{D}$$ is a
unitary isomorphism that intertwines the $G$-actions.
\item  In the case $\mathbb{H}=\bH_\nu$ we have $\bV=\mathbb{C}1$,
where $1$ is the constant function, and $$(q_T 1, F)=(1,
F(T))=\overline{F(T)}, \quad T\in \mathcal{D}.$$
\end{enumerate}
\end{theorem}

In the case where $\mathbb{H}= L^2(\Omega,d\mu_\nu)$ the map $\Xi$
is given in Lemma \ref{le-uniiso}:
\begin{eqnarray*}
\Xi (f)(w)&=&\frac{1}{\sqrt{\Gamma_\gO (\nu)}}\, \pi_\nu (c)^{-1}\cL (f)(w)\\
&=&\sqrt{\frac{2^{r\nu}}{\Gamma_\gO (\nu)}}\, \int_\gO e^{-(c
w,x)}f(x)\Delta(x)^{\nu-d/r}dx
\end{eqnarray*} The highest weight space is
$H=\mathbb{C}\ell_\circ^\nu$. Let $v_\circ=\sqrt
{\frac{2^{r\nu}}{\Gamma_{\Omega}(\nu)}}\ell_\circ^\nu$. Then
$v_\circ$ has norm $1$ and $\Xi v_\circ$ is the constant function
$1$ on $\mathcal{D}$.

\begin{theorem}[\cite{df2}] Let $\bE$ be a closed subspace of $L^2(\Omega,d\mu_\nu)$,
$\set{e_{{ \alpha}}}$ an orthonormal basis of $\bE$, and
$E_{ \alpha}=\Xi(e_{ \alpha})\in \bH_\nu$ for each
$\alpha$. Let $\mathrm{pr}_{\bE}$ be the orthogonal projection of
$L^2(\Omega,d\mu_\nu)$ onto $\bE$.  Then we have for $w\in
\mathcal{D}$ and $x\in \Omega$:
$$ \beta_\nu \,\overline{\Delta(e-w)^{-\nu}}\mathrm{pr}_{\bE}(e^{-
\overline{(\bc(w), x)}})=\sum_{{
\alpha}}\overline{E_{\alpha}(w)}e_{ \alpha}(x),$$ where
$\beta_\nu=\sqrt{\frac{2^{r\nu}}{\Gamma_\gO (\nu)}}$ and the
convergence is with respect to the Hilbert space norm in
$L^2(\Omega,d\mu_\nu).$
\end{theorem}
\begin{proof}
First we have
\begin{eqnarray*}
\mathrm{pr}_{\bE}(q_w(v_\circ)) &=&\sum_{\alpha}(q_w(v_\circ), e_\alpha)e_\alpha\\
&=&\sum_{\alpha}(\Xi q_w v_\circ, \Xi e_\alpha)e_\alpha\\
&=&\sum_{\alpha}(q_w\Xi v_\circ,  E_\alpha)e_\alpha\\
&=& \sum_{\alpha }(1 ,  E_\alpha(w))e_\alpha\\
&=&\sum_{\alpha}\overline{E_\alpha(w)}e_\alpha
\end{eqnarray*}
On the other hand,
\begin{eqnarray*}
(q_w v_\circ, f)&=&(\Xi q_w v_\circ, \Xi f)\\
&=& (q_w 1, \Xi f)\\
&=& \overline {\Xi f(w)}\\
&=&\beta_\nu \;\overline {\Delta(e-w)^{-\nu}\int_\Omega e^{-(c w,
x)}f(x)
\, d\mu_\nu(x)}\\
&=& \beta_\nu\;\overline{\Delta(e-w)^{-\nu}}\overline{(e^{-(c
w,x)}, f)}\
\end{eqnarray*}
From this it follows that

$$q_w v_\circ(x)=\beta_\nu\;\overline{\Delta(e-w)^{-\nu}} e^{\overline{-(\bc (w) ,x)}}.$$
and
$$\beta_\nu\;\overline{\Delta(e-w)^{-\nu}}\mathrm{pr}_{\bE}(e^{-
\overline{(\bc (w) , x)}})=\sum_{{
\alpha}}\overline{E_{\alpha}(w)}e_{\alpha}(x)\, .$$
\end{proof}

We now specialize to the case where
$\bE=L^2(\Omega,d\mu_\nu)^{K\cap H}$ with orthonormal basis
$$e_{m}= \sqrt{\frac{2^{r\nu}d_{m}}{\Gamma_\gO ( \nu)
\left(\frac{d}{r}\right)_{m}(\nu )_{m}}}\,
\ell^\nu_{m}.
$$ The orthogonal projection
$\mathrm{pr}_{\bE}:L^2(\Omega,d\mu_\nu)\rightarrow \bE$ is given
by $$\mathrm{pr}_{\bE}(f)=\int_{K\cap H}f(kx)\;dk.$$   The above
theorem immediately gives:

\begin{theorem}[\cite{df2}]
Let $w \in \mathcal{D}$ and $x\in \Omega$. Then
$$ \Delta(e-w)^{-\nu}\int_{K\cap H}
e^{-(k\cdot x,  (1+w)(1-w)^{-1})}\; dk=\sum_{{m}\in
\Lambda}d_{m}\frac{1}{(\frac n r)_{{ m}}}
\psi_{m}(w)\ell_{m}^\nu(x).$$
\end{theorem}

This formula is the generating function for the Laguerre functions
on symmetric cones and generalizes equation \eqref{eq-generating2}
(cf \cite{FK} p.347 and the references given there).

\providecommand{\bysame}{\leavevmode\hbox
to3em{\hrulefill}\thinspace}

\end{document}